\documentclass[12pt]{elsarticle}

\usepackage{amsmath,amsthm,amssymb,a0size,bm,indentfirst,txfonts}
\usepackage{stmaryrd,CJK,color,graphicx,palatino}
\usepackage[toc,page,title,titletoc,header]{appendix}
\usepackage{paralist}
\usepackage{epsfig}
\usepackage[colorlinks=true]{hyperref}
\hypersetup{citecolor=red}

\usepackage{remreset}

\allowdisplaybreaks

\newcommand{\bex}{\begin{eqnarray*}}
\newcommand{\eex}{\end{eqnarray*}}
\newcommand{\be}{\begin{eqnarray}}
\newcommand{\ee}{\end{eqnarray}}
\newcommand{\ba}{\begin{array}}
\newcommand{\ea}{\end{array}}
\newcommand{\bi}{\begin{itemize}}
\newcommand{\ei}{\end{itemize}}
\newcommand{\bn}{\begin{enumerate}}
\newcommand{\en}{\end{enumerate}}

\newcommand{\sex}[1]{\left(#1\right)}

\newcommand{\sez}[1]{\left[#1\right]}
\newcommand{\sed}[1]{\left\{#1\right\}}
\newcommand{\sev}[1]{\left|#1\right|}
\newcommand{\sen}[1]{\left\Vert#1\right\Vert}

\newtheorem{thm}{\indent Theorem}

\newtheorem{lem}[thm]{\indent Lemma}

\newtheorem{defn}[thm]{\indent Definition}
\newtheorem{rem}[thm]{\indent Remark}

\newcommand{\bp}{\begin{proof}}
\newcommand{\ep}{\end{proof}}

\newcommand{\R}{\bm{R}}

\def\nGamma{\mathnormal{\Gamma}}
\def\ds{\mathrm{d}s}
\def\dx{\mathrm{d}x}
\def\n{\nabla}
\def\p{\partial}
\def\lap{\Delta}

\def\bbb{\bm{b}}

\def\bbu{\bm{u}}

\def\bbomega{\bm{\omega}}

\def\curl{\mbox{curl }}

\def\ve{\varepsilon}

\makeatletter
\def\@eqnnum{{\normalfont \color{red} (\theequation)}}
\makeatother

\makeindex

\begin{document}

\begin{frontmatter}

\title{A regularity criterion for the Navier-Stokes equations via two entries of the velocity Hessian tensor}

\author{Zujin Zhang\fnref{z}\corref{cor}}

\address{
Department of Mathematics,
Sun Yat-sen University\\
Guangzhou 510275, Guangdong, P.R. China}

\cortext[cor]{Corresponding author}

\fntext[z]{uia.china@gmail.com}

\begin{abstract}
We consider the Cauchy problem for the incompressible Navier-Stokes equations in $\R^3$, and provide a sufficient condition to ensure the smoothness of the solution. It involves only two entries of the velocity Hessian.
\end{abstract}

\begin{keyword}
Incompressible Navier-Stokes equations\sep
regularity criterion\sep
global regularity\sep
weak solutions\sep
strong solutions

\MSC[2010]
35Q30\sep
35B45\sep
76D03

\end{keyword}

\end{frontmatter}

\section{Introduction}
\label{sect:intro}
This paper is concerned with the global regularity of solutions of the three-dimensional Navier-Stokes equation (NSE):
\be\label{1:eq:NSE}
\left\{
\ba{lll}
\p_t\bbu-\nu\lap \bbu
+\sex{\bbu\cdot\n}\bbu+\n p=0,
&\mbox{in }\R^3\times(0,T),\\
\n\cdot\bbu=0,
&\mbox{in }\R^3\times (0,T),\\
\bbu=\bbu_0,
&\mbox{on }\R^3\times \sed{t=0},
\ea
\right.
\ee
where $T>0$ is a given time, $\bbu=\sex{u_1,u_2,u_3}$ is the velocity field, $p$ is a scalar pressure, and $\bbu_0$ is the initial velocity field satisfying $\n\cdot \bbu_0=0$ in the sense of distributions.

The global existence of a weak solution $\bbu$ to \eqref{1:eq:NSE} with initial data of finite energy is well-known since the work of Leary \cite{Leray_34}, see also Hopf \cite{Hopf}. However, the issue of uniqueness and regularity of $\bbu$ was left open, and is still unsolved up to date. Pioneered by Serrin \cite{Serrin_62,Serrin_63} and Prodi \cite{Prodi_59}, there have been a lot of literatures devoted to finding sufficient conditions to ensure the smoothness of $\bbu$. These conditions involve either
\bi
\item[--]the velocity $\bbu$, see \cite{Escauriaza_Seregin_Sverak_03,FJR_72,Giga_86,
    Serrin_62,Serrin_63,Sohr_83,Sohr_Wahl_86}, which states
\be\label{1:Serrin}
\bbu\in L^\alpha(0,T;L^\beta(\R^3)),\ \ \
\indent
\frac{2}{\alpha}+\frac{3}{\beta}=1,\ \ \
\indent
3\leq \beta\leq \infty;
\ee
\item[--]or several components of the velocity $\bbu$, the velocity gradient $\n\bbu$, the vorticity $\bbomega=\curl \bbu$, or the pressure gradient $\n p$, see
\cite{Beirao da Veiga_95,
Beirao da Veiga_02,
BerselliGaldi_02,
Cao_2010_one component,
Cao_Titi_08_one component,
Cao_Titi_10_one entry reg,
Chae_Lee_01,
Constantin_Fefferman_93,
Escauriaza_Seregin_Sverak_03,
Fan_Jiang_Ni_08,
Kim_10,
Kukavica_Ziane_07_one direction,
Kukavica_Ziane_06_one component,
Neustupa_Novotny_Penel_99,
Neustupa_Penel,
Penel_Pokorny_10,
Penel_Pokorny,
Zhou_05_direction vorticity,
Zhou_02_grad u3,
Zhou_05_one component,
Zhou_05_direction_ANZIAM,
Zhou_06_grad pressure,
Zhou_06_pressure,
Zhou_04_pressure,
Zhang_Xicheng_08,
Zhang_Chen_05,
Zhang_10_some,Zhang_10_two}, and references cited therein.
\ei

We remark that many of the regularity criteria established in the above cited papers have been extended to the following three dimensional magneto-hydrodynamic equations (MHD):
\be\label{1:MHD}
\left\{\ba{llll}
\p_t\bbu+\sex{\bbu\cdot\n}\bbu
-\sex{\bbb\cdot\n}\bbb
-\lap\bbu+\n p=0,
&\mbox{in }\R^3\times(0,\infty),\\
\p_t\bbu+\sex{\bbu\cdot\n}\bbb
-\sex{\bbb\cdot\n}\bbu
-\lap \bbb=0,
&\mbox{in }\R^3\times(0,\infty),\\
\n\cdot\bbu=0,\ \n\cdot\bbb=0,
&\mbox{in }\R^3\times(0,\infty),\\
\bbu=\bbu_0,\ \bbb=\bbb_0,
&\mbox{on }\R^3\times\sed{t=0}.
\ea\right.
\ee
Here $\bbu=(u_1,u_2,u_3)$ is the velocity field, $\bbb=(b_1,b_2,b_3)$ is the magnetic field, $\bbu_0$, $\bbb_0$ are the corresponding initial data, and $p$ is a scalar pressure. The interested readers are referred to \cite{HeXin_05,WuJH_03,WuJH_08,Zhou_06_MHD pressure,Zhou_05_MHD_grad u,Zhang_10_JMAA} and references cited therein.

Motivated by \cite{Cao_Wu_2010_MHD} and \cite{Zhou_Pokorny_09_one component}, we consider in this paper the regularity criterion involving $\p_1\p_3u_3$ and $\p_2\p_3u_3$ only. Before stating the precise result, let us recall the weak formulation of \eqref{1:eq:NSE}.
\begin{defn}
Let $\bbu_0\in L^2(\R^3)$ with $\n\cdot u_0=0$, and $T>0$. A measurable $\R^3$-valued vector $\bbu$ is said to be a weak solution of \eqref{1:eq:NSE} if the following conditions hold:
\bn
\item $\bbu\in L^\infty(0,T;L^2(\R^3))\cap L^2(0,T;H^1(\R^3))$;
\item $\bbu$ solves $\eqref{1:eq:NSE}_{1,2}$ in the sense of distributions; and
\item the energy inequality, that is,
    \be\label{1:eq:energy}
    \sen{u(t)}_2^2
    +\nu\int_{t_0}^t\sen{\n\bbu(s)}_2^2ds
    \leq \sen{\bbu(t_0)}_2^2,
    \ee
    for almost every $t_0$ (including $t_0=0$) and every $t\geq t_0$.
\en
\end{defn}

Our regularity criterion now reads:
\begin{thm}\label{1:thm:main}
Let $\bbu_0\in H^1(\R^3)$ with $\n\cdot\bbu_0=0$, and $T>0$. Suppose $\bbu$ is the corresponding solution on $[0,T]$ of \eqref{1:eq:NSE}. If
\be\label{1:thm:main:reg}
\p_1\p_3u_3,\ \p_2\p_3u_3\ \in
L^\alpha(0,T;L^\beta(\R^3)),\indent \frac{2}{\alpha}+\frac{3}{\beta}=2+\frac{1}{\beta},
\indent 1< \beta\leq\infty,
\ee
then $\bbu$ is smooth on $(0,T)$.
\end{thm}

Before proving this theorem in Section \ref{2}, we collect here some notations used throughout this paper and make some remarks on our result.

The usual Lebesgue spaces $L^q(\R^3)$ $(1\leq q\leq \infty)$ is endowed with the norm $\sen{\cdot}_{q}$.
For a Banach space $\sex{X,\sen{\cdot}}$, we do not distinguish it with its vector analogues $X^3$, thus the norm in $X^3$ is still denoted by $\sen{\cdot}$; however, all vector- and tensor-valued functions are printed boldfaced. We also denote by
\bex
\p_i \varphi=\frac{\p \varphi}{\p x_i},\indent
\p_{ij}^2\varphi=\p_i\p_j\varphi
=\frac{\p^2 \varphi}{\p x_i\p x_j}
,\indent 1\leq i\leq 3
\eex
the first- and second-order derivatives of a function $\varphi$; by
\bex
\n_h \varphi=\sex{\p_1,\p_2}\varphi,\ \lap_h \varphi=\sex{\p_{11}^2+\p_{22}^2}\varphi
\eex
the horizontal gradient, horizontal Laplacian of $\varphi$.

\begin{rem}
Noticing that
\bex
\lim_{\beta\to 1^-}\sex{2+\frac{1}{\beta}}=3,
\eex
we almost establish a Serrin-type regularity criterion via $\p_1\p_3u_3$ and $\p_2\p_3u_3$ only.
\end{rem}

\begin{rem}
Due to the orthogonal transformation invariance of NSE \eqref{1:eq:NSE}, we easily extends our regularity criterion as
\bex
\p_i\p_ku_k,\ \p_j\p_ku_k\ \in
L^\alpha(0,T;L^\beta(\R^3)),\indent \frac{2}{\alpha}+\frac{3}{\beta}=2+\frac{1}{\beta},\indent
1< \beta\leq\infty,
\eex
where $\sed{i,j,k}=\sed{1,2,3}$.
\end{rem}
\begin{rem}
Our result seems to be more involved than that in \cite{Cao_Titi_10_one entry reg} in the following sense. In \cite{Cao_Titi_10_one entry reg}, only conditions on $1$ component is needed among the total $9$ components of the velocity gradient tensor
\bex
\p\otimes \bbu=\n\bbu=\sez{\p_iu_j},
\eex
the ratio is $1/9$. And our result requires regularity of $2$ entries  among the total $27$ entries of the velocity Hessian tensor
\bex
\p^2\otimes \bbu=\n^2\bbu=\sez{\p_{ij}^2u_k},
\eex
the ratio being $2/27$, which is less than $1/9$.
\end{rem}
\begin{rem}
In classical and numerical analysis, we do not only rely on the graph of the differential to study properties of a function, but also utilize the graph of the Hessian to do so. In this point of view, our result is a complement to that in \cite{Cao_Titi_10_one entry reg,Zhou_Pokorny_09_one component}.
\end{rem}
\begin{rem}
Our proof in Section \ref{2} is different from that in \cite{Cao_Titi_10_one entry reg,Zhou_Pokorny_09_one component}, due to the our assumptions on the Hessian. We shall first bound $\sen{\n_h\bbu}_2$, then estimate $\sen{\n\bbu}_2$. In fact, tracking the proof of that in \cite{Cao_Titi_10_one entry reg,Zhou_Pokorny_09_one component} we will obtain a not-so-good result.
\end{rem}

\section{Proof of the main result}\label{2}
In this section, we shall prove Theorem \ref{1:thm:main}. First, let us recall and prove some technical lemmas.

The first one being a component-reducing technique due to Kukavica and Ziane \cite{Kukavica_Ziane_07_one direction}.
\begin{lem}\label{2:lem:kukavica_ziane}
Assume $\bbu=\sex{u_1,u_2,u_3}\in C_c^\infty(\R^3)$ is divergence free. Then
\be\label{lem:kukavica_ziane:eq:u_3}
\sum_{i,j=1}^2\int_{\R^3}u_i\p_iu_j\lap_hu_j\dx
&=&\frac{1}{2}\sum_{i,j=1}^2
\int_{\R^3}\p_iu_j\p_iu_j\p_3u_3\dx\nonumber\\
& &
-\int_{\R^3}\p_1u_1\p_2u_2\p_3u_3\dx
+\int_{\R^3}\p_1u_2\p_2u_1\p_3u_3\dx.
\ee
\end{lem}
And the next two lemmas are variants of multiplicative Sobolev inequalities in $\R^3$, similar in spirit to that in \cite{Cao_Titi_10_one entry reg}.
\begin{lem}\label{2:lem:3}
Let $1<r\leq 3$. Assume $f,\ g,\ h\ \in C_c^\infty(\R^3)$. Then there exists a constant $C>0$ such that
\bex
\int_{\R^3}f\ g\ h\ \dx
\leq C \sen{f}_2^\frac{r-1}{r}
\sen{\p_3 f}_\frac{2}{3-r}^\frac{1}{r}
\sen{g}_2^\frac{r-1}{r}
\sen{\p_1 g}_2^\frac{1}{2r}
\sen{\p_2 g}_2^\frac{1}{2r}
\sen{h}_2^\frac{r-1}{r}
\sen{\p_1 h}_2^\frac{1}{2r}
\sen{\p_2 h}_2^\frac{1}{2r}.
\eex
\end{lem}
\bp
\bex
& &\int_{\R^3}f\ g\ h\ \dx_1 \dx_2 \dx_3\\
& &\leq \int_{\R^2}
\max_{x_3}\sev{f}
\cdot
\sex{\int_{\R} \sev{g}^2\dx_3}^{1/2}
\cdot
\sex{\int_{\R} \sev{h}^2\dx_3}^{1/2}
\dx_1\dx_2\\
& &
\leq
\sez{\int_{\R^2}
\sex{\max_{x_3}\sev{f}}^r
\dx_1\dx_2}^{1/r}
\cdot
\sez{
\int_{\R^2}
\sex{\int_{\R}\sev{g}^2dx_3}^\frac{r}{r-1}
\dx_1\dx_2
}^\frac{r-1}{2r}\\
& &\ \ \cdot
\sez{
\int_{\R^2}
\sex{\int_{\R}\sev{h}^2dx_3}^\frac{r}{r-1}
\dx_1\dx_2
}^\frac{r-1}{2r}\\
& &\leq C
\sez{
\int_{\R^3}
\sev{f}^{r-1}
\sev{\p_3 f}\dx
}^{1/r}
\cdot
\sez{
\int_{\R}
\sex{
\int_{\R^2}\sev{g}^\frac{2r}{r-1}\dx_1\dx_2
}^\frac{r-1}{r}
\dx_3
}^{1/2}\\
& &\ \ \cdot
\sez{
\int_{\R}
\sex{
\int_{\R^2}\sev{h}^\frac{2r}{r-1}\dx_1\dx_2
}^\frac{r-1}{r}
\dx_3
}^{1/2}\\
& &
\leq C \sen{f}_2^\frac{r-1}{r}
\sen{\p_3 f}_\frac{2}{3-r}^\frac{1}{r}
\sen{g}_2^\frac{r-1}{r}
\sen{\p_1 g}_2^\frac{1}{2r}
\sen{\p_2 g}_2^\frac{1}{2r}
\sen{h}_2^\frac{r-1}{r}
\sen{\p_1 h}_2^\frac{1}{2r}
\sen{\p_2 h}_2^\frac{1}{2r}.
\eex
\ep
The same argument also yields
\begin{lem}\label{2:lem:1 or 2}
Let $1<r\leq 3$. Assume $f,\ g,\ h\ \in C_c^\infty(\R^3)$. Then there exists a constant $C>0$ such that
\bex
\int_{\R^3}f\ g\ h\ \dx
\leq C \sen{f}_2^\frac{r-1}{r}
\sen{\p_1 f}_\frac{2}{3-r}^\frac{1}{r}
\sen{g}_2^\frac{r-1}{r}
\sen{\p_2 g}_2^\frac{1}{2r}
\sen{\p_3 g}_2^\frac{1}{2r}
\sen{h}_2^\frac{r-1}{r}
\sen{\p_2 h}_2^\frac{1}{2r}
\sen{\p_3 h}_2^\frac{1}{2r}.
\eex
\end{lem}
\bigskip
{\bf Proof of Theorem \ref{1:thm:main}.}

\bigskip
{\bf Step 1.} Preliminary reduction.
\bigskip

For any $\ve\in (0,T)$, due to the fact that $\n\bbu\in L^2(0,T;L^2(\R^3))$, we may find a $\delta\in (0,\ve)$, such that $\n\bbu(\delta)\in L^2(\R^3)$. Take this $\bbu(\delta)$ as initial data, there exists an $\tilde\bbu\in C([\delta,\nGamma^*),H^1(\R^3)) \cap L^2(0,\nGamma^*;H^2(\R^3))$, where $[\delta, \nGamma^*)$ is the life span of the unique strong solution, see \cite{Temam_77_NSE}. Moreover, $\tilde\bbu\in C^\infty(\R^3\times (\delta,\nGamma^*))$. According to the uniqueness result, $\tilde \bbu=\bbu$ on $[\delta,\nGamma^*)$. If $\nGamma^*\geq T$, we have already that $\bbu\in C^\infty(\R^3\times (0,T))$, due to the arbitrariness of $\ve\in (0,T)$. In case $\nGamma^*<T$, our strategy is to show that $\sen{\n\bbu(t)}_2$ remains bounded independently of $t\nearrow\nGamma^*$. The standard continuation argument then yields that $[\delta,\nGamma^*)$ could not be the maximal interval of existence of $\tilde\bbu$, and consequently $\nGamma^*\geq T$. This concludes the proof.

For this purpose, let us choose a $\tau$ sufficiently close to $\nGamma^*$ such that
\be\label{2:small}
\int_\tau^{\nGamma^*} \sez{
\sen{\sex{\p_1,\p_2}\p_3u_3(s)}_\beta^\alpha
+\sen{\n\bbu(s)}_2^2
}ds<\tilde \ve,
\ee
where $\tilde \ve>0$ is small and will be chosen later on.

We shall first in Step $2$ establish the bounds of $\sen{\n_h\bbu(t)}_2$ for $t\in [\tau,\nGamma^*)$. The estimation of $\sen{\n\bbu(t)}_2$ is derived in Step $3$.

\bigskip
{\bf Step 2.} $\sen{\n_h\bbu(t)}_2$ estimates.
\bigskip

Taking the inner product of $\eqref{1:eq:NSE}_1$ with $-\lap_h\bbu$, we have
\be
\label{2:I}
& &\frac{1}{2}\frac{d}{dt}
\sen{\n_h\bbu}_2^2
+\nu\sen{\n\n_h\bbu}_2^2
=\int_{\R^3}\sex{\bbu\cdot\n}\bbu\cdot\lap_h\bbu \dx
\nonumber\\
& &=\sum_{i,j=1}^2\int_{\R^3}u_i\p_i\lap_hu_j\dx
+\sum_{j=1}^2\int_{\R^3}u_3\p_3u_j\lap_hu_j\dx
+\sum_{i=1}^3\int_{\R^3}u_i\p_iu_3\lap_hu_3\dx\nonumber\\
& &\equiv I_1+I_2+I_3.
\ee
Invoking Lemma \ref{2:lem:kukavica_ziane}, $I_1$ can be rewritten as
\be\label{2:I_1}
I_1=\frac{1}{2}\sum_{i,j=1}^2
\int_{\R^3}\p_iu_j\p_iu_j\p_3u_3\dx
-\int_{\R^3}\p_1u_1\p_2u_2\p_3u_3\dx
+\int_{\R^3}\p_1u_2\p_2u_1\p_3u_3\dx.
\ee
For $I_2,\ I_3$, it follows by integrating by parts and noticing the divergence free condition $\n\bbu=0$, that
\be\label{2:I_2}
I_2&=&\sum_{j=1}^2\int_{\R^3}u_3\p_3u_j\lap_hu_j\dx\nonumber\\
&=&-\sum_{j,k=1}^2\int_{\R^3}\p_ku_3\p_3u_j\p_ku_j\dx
+\frac{1}{2}\sum_{j,k=1}^2\int_{\R^3}\p_3u_3\sex{\p_ku_j}^2\dx,
\ee
\be\label{2:I_3}
I_3&=&\sum_{i=1}^3\int_{\R^3}u_i\p_iu_3\lap_h u_3\dx\nonumber\\
&=&\sum_{i=1}^3\sum_{j=1}^2
\int_{\R^3}\p_ju_i\p_iu_3\p_ju_3\dx.
\ee
Gathering \eqref{2:I_1},\eqref{2:I_2},\eqref{2:I_3} together, \eqref{2:I} becomes
\be\label{2:J}
\frac{1}{2}\frac{d}{dt}
\sen{\n_h\bbu}_2^2
+\nu\sen{\n\n_h\bbu}_2^2
&\leq&
C\int_{\R^3}\sev{\n_hu_3}\cdot\sev{\n\bbu}\cdot\sev{\n_h\bbu}\dx
+C\int_{\R^3}\sev{\p_3u_3}\cdot\sev{\n_h\bbu}^2\dx\nonumber\\
&\equiv&J_1+J_2.
\ee
We now apply Lemmas \ref{2:lem:3} and \ref{2:lem:1 or 2} with $\beta=\frac{2}{3-r}$ to bound $J_1$, $J_2$ respectively as
\be\label{2:J_1}
J_1&=&C\int_{\R^3}\sev{\n_hu_3}\cdot\sev{\n\bbu}\cdot\sev{\n_h\bbu}\dx
\nonumber\\
&\leq&C\sen{\n_hu_3}_2^\frac{2(\beta-1)}{3\beta-2}
\sen{\n_h\p_3u_3}_\beta^\frac{\beta}{3\beta-2}
\sen{\n\bbu}_2^\frac{2(\beta-1)}{3\beta-2}
\sen{\n_h\bbu}_2^\frac{2(\beta-1)}{3\beta-2}
\sen{\n\n_h\bbu}_2^\frac{2\beta}{3\beta-2},
\ee
\be\label{2:J_2}
J_2& =&
C\int_{\R^3}\sev{\p_3u_3}\cdot\sev{\n_h\bbu}^2\dx\nonumber\\
&\leq&C\sen{\p_3u_3}_2^\frac{2(\beta-1)}{3\beta-2}
\sen{\p_1\p_3u_3}_\beta^\frac{\beta}{3\beta-2}
\sen{\n_h\bbu}_2^\frac{4(\beta-1)}{3\beta-2}
\sen{\n\n_h\bbu}_2^\frac{2\beta}{3\beta-2}\nonumber\\
&\leq&C\sen{\n_h\bbu}_2^\frac{6(\beta-1)}{3\beta-2}
\sen{\p_1\p_3u_3}_\beta^\frac{\beta}{3\beta-2}
\sen{\n\n_h\bbu}_2^\frac{2\beta}{3\beta-2},
\ee
where in the last inequality, we use the divergence free condition $\n\cdot\bbu=0$.

Substituting \eqref{2:J_1}, \eqref{2:J_2} into \eqref{2:J}, we obtain by Young inequality that
\bex
\frac{1}{2}\frac{d}{dt}
\sen{\n_h\bbu}_2^2
+\nu \sen{\n\n_h\bbu}_2^2
&\leq& C\sen{\n_h\bbu}_2^\frac{4(\beta-1)}{3\beta-2}
\sen{\n\bbu}_2^\frac{2(\beta-1)}{3\beta-2}
\sen{\sex{\p_1,\p_2}\p_3u_3}_\beta^\frac{\beta}{3\beta-2}
\sen{\n\n_h\bbu}_2^\frac{2\beta}{3\beta-2}\\
&\leq&
C\sen{\n_h\bbu}_2^2
\sen{\n\bbu}_2
\sen{\sex{\p_1,\p_2}\p_3u_3}_\beta^\frac{\beta}{2(\beta-1)}
+\frac{\nu}{2}\sen{\n\n_h\bbu}_2^2,
\eex
i.e.
\be
\label{2:grad_h u_diff}
\frac{d}{dt}
\sen{\n_h\bbu}_2^2
+\nu \sen{\n\n_h\bbu}_2^2
\leq C\sen{\n_h\bbu}_2^2
\sen{\n\bbu}_2
\sen{\sex{\p_1,\p_2}\p_3u_3}_\beta^\frac{\beta}{2(\beta-1)}.
\ee
Integrating this inequality over $[\tau,t]$, for any $t\in [\tau,\nGamma^*)$, we get
\be
\label{2:grad_h u:integ}
& &\sen{\n_h\bbu(t)}_2^2
+\nu\int_\tau^t\sen{\n\n_h\bbu(s)}_2^2\ds\nonumber\\
& &\leq
\sen{\n_h\bbu(\tau)}_2^2
+C\int_\tau^t
\sen{\sex{\p_1,\p_2}\p_3u_3(s)}_\beta
^\frac{\beta}{2(\beta-1)}
\sen{\n_h\bbu(s)}_2^2
\sen{\n\bbu(s)}_2\ds.
\ee
Young inequality then implies
\bex
\sup_{\tau\leq t<\nGamma^*}\sen{\n_h\bbu(t)}_2^2
&\leq&
\sen{\n_h\bbu(\tau)}_2^2\\
& &
+C
\sup_{\tau\leq t<\nGamma^*}
\sen{\n_h\bbu(t)}_2^2
\cdot \int_\tau^{\nGamma^*}
\sez{\sen{\sex{\p_1,\p_2}\p_3u_3(s)}
_\beta^\frac{\beta}{\beta-1}
+\sen{\n\bbu(s)}_2^2}ds.
\eex
Hence if $\tilde \ve$ in \eqref{2:small} is choose so that
\bex
C\tilde\ve<\frac{1}{2},
\eex
we see
\be\label{2:grad_h u:bdd}
\sup_{\tau\leq t<\nGamma^*}
\sen{\n_h\bbu(t)}_2^2
\leq 2\sen{\n_h\bbu(\tau)}_2^2.
\ee

\bigskip
{\bf Step 3.} $\sen{\n\bbu(t)}_2$ estimates.
\bigskip

Taking the inner product of $\eqref{1:eq:NSE}_1$ with $-\lap\bbu$ in $L^2(\R^3)$, we see
\be\label{2:K}
\frac{1}{2}\frac{d}{dt}
\sen{\n\bbu}_2^2
+\nu\sen{\lap\bbu}_2^2
&=&\int_{\R^3}\sex{\bbu\cdot\n\bbu}\cdot\lap\bbu \dx
\nonumber\\
&=&\int_{\R^3}\sex{\bbu\cdot\n\bbu}\cdot\lap_h\bbu \dx
+\int_{\R^3}\sex{\bbu\cdot\n\bbu}\cdot\p_{33}^2\bbu \dx\nonumber\\
&\equiv&K_1+K_2.
\ee
The term $K_1$ can be dominated similarly as in Step $2$, and we find
\be\label{2:K_1}
K_1&=&\int_{\R^3}\sex{\bbu\cdot\n\bbu}\cdot\lap_h\bbu \dx\nonumber\\
&\leq&C\sen{\sex{\p_1,\p_2}\p_3u_3}
_\beta^\frac{\beta}{2(\beta-1)}
\sen{\n_h\bbu}_2^2\sen{\n\bbu}_2
+\frac{\nu}{4}\sen{\lap\bbu}_2^2.
\ee
Meanwhile for $K_2$, we have, by integrating by parts and noticing the divergence free condition $\n\cdot\bbu=0$, that
\bex
K_2&=&\int_{\R^3}\sex{\bbu\cdot\n}\bbu\cdot\p_{33}^2\bbu\dx\\
&=&\sum_{i,j=1}^3\int_{\R^3}\p_3u_i\p_iu_j\p_3u_j\dx\\
&=&\sum_{i=1}^2\sum_{j=1}^3
\int_{\R^3}\p_3u_i\p_iu_j\p_3u_j\dx
+\sum_{j=1}^3\int_{\R^3}\p_3u_3\p_3u_j\p_3u_j\dx\\
&=&\sum_{i=1}^2\sum_{j=1}^3\int_{\R^3}\p_3u_i\p_iu_j\p_3u_j\dx
-\sum_{j=1}^3\sex{\p_1u_1+\p_2u_2}\p_3u_j\p_3u_j\dx.
\eex
Applying H\"older inequality, interpolation inequality, Sobolev inequality and Young inequality yields
\be
\label{2:K_2}
K_2&\leq&C\int_{\R^3}\sev{\n_h\bbu}\cdot\sev{\n\bbu}^2\dx\nonumber\\
&\leq&C\sen{\n_h\bbu}_2\sen{\n\bbu}_4^2\nonumber\\
&\leq&C\sen{\n_h\bbu}_2\sen{\n\bbu}_2^{1/2}
\sen{\lap\bbu}_2^{3/2}\nonumber\\
&\leq&C\sen{\n_h\bbu}_2^4
\sen{\n\bbu}_2^2
+\frac{\nu}{4}\sen{\lap\bbu}_2^2.
\ee
Now, replacing \eqref{2:K_1}, \eqref{2:K_2} into \eqref{2:K}, we get
\bex
\frac{d}{dt}\sen{\n\bbu}_2^2
+\nu\sen{\lap\bbu}_2^2
\leq C\sen{\sex{\p_1,\p_2}\p_3u_3}_\beta^\frac{\beta}{2(\beta-1)}
\sen{\n_h\bbu}_2^2\sen{\n\bbu}_2
+C\sen{\n_h\bbu}_2^4
\sen{\n\bbu}_2^2.
\eex
Thanks to \eqref{2:grad_h u:bdd}, we obtain further that
\bex
\frac{d}{dt}\sen{\n\bbu}_2^2
+\nu\sen{\lap\bbu}_2^2
&\leq& C\sen{\n_h\bbu}^2\sez{1+
\sen{\sex{\p_1,\p_2}\p_3u_3}_\beta^\frac{\beta}{\beta-1}
\sen{\n\bbu}_2^2
}\\
& &
+4C\sen{\n_h\bbu(\tau)}_2^4\sen{\n\bbu}_2^2,\quad \mbox{on }[\tau,\nGamma^*).
\eex
Invoking Gronwall inequality then implies $\sen{\n\bbu(t)}_2$, $t\in [\tau,\nGamma^*)$ is uniformly bounded, as desired. This completes the proof of Theorem \ref{1:thm:main}.
\hfill $\square$


\begin{thebibliography}{00}

\bibitem{Beirao da Veiga_95}
H. Beir\~ao da Veiga,
\emph{A new regularity class for the Navier-Stokes equations in $\R^n$,}
Chinese Ann. Math. Ser. B
{\bf 16}
(1995),
407--412.

\bibitem{Beirao da Veiga_02}
H. Beir\~ao da Veiga, L.C. Berselli,
\emph{On the regularizing effect of the vorticity direction in incompressible viscous flows,}
Differential Integral Equations
{\bf 15}
(2002),
345--356.

\bibitem{BerselliGaldi_02}
L.C. Berselli, G.P. Galdi,
\emph{Regularity criteria involving the pressure for the weak solutions to the Navier-Stokes equations,}
Proc. Amer. Math. Soc.
{\bf 130}
(2002),
3585--3595.

\bibitem{Cao_2010_one component}
C.S. Cao,
\emph{Sufficient conditions for the regularity to the $3$D Navier-Stokes equations,}
Discrete Contin. Dyn. Syst.
{\bf 26}
(2010),
1141--1151.

\bibitem{Cao_Titi_10_one entry reg}
C.S. Cao, E.S. Titi,
\emph{Global regularity criterion for the $3$D Navier-Stokes equations involving one entry of the velocity gradient tensor,}
arXiv: 1005.4463 [math. AP] 25 May 2010.

\bibitem{Cao_Titi_08_one component}
C.S. Cao, E.S. Titi,
\emph{Regularity criteria for the three-dimensional Navier-Stokes equations,}
Indiana Univ. Math. J.
{\bf 57}
(2008),
2643--2661.

\bibitem{Cao_Wu_2010_MHD}
C.S. Cao, J.H. Wu,
\emph{Two new regularity criteria for the $3$D MHD equations,}
J. Differential Equations
{\bf 248}
(2010),
2263--2274.

\bibitem{Chae_Lee_01}
D. Chae, J. Lee,
\emph{Regularity criterion in terms of pressure for the Navier-Stokes equations,}
Nonlinear Anal. TMA
{\bf 46}
(2001),
727--735.

\bibitem{Constantin_Fefferman_93}
P. Constantin, C. Fefferman,
\emph{Direction of vorticity and the problem of global regularity for the Navier-Stokes equations,}
Indiana Univ. Math. J.
{\bf 42}
(1993),
775--789.

\bibitem{Escauriaza_Seregin_Sverak_03}
L. Escauriaza, G. Seregin, V. \v Sver\'ak,
\emph{Backward uniqueness for parabolic equations,}
Arch. Ration. Mech. Anal.
{\bf 169}
(2003),
147--157.

\bibitem{FJR_72}
E.B. Fabes,B.F. Jones, N.M. Rivi\'ere,
\emph{The initial value problem for the Navier-Stokes equations with data in $L^{p}$,}
Arch. Rational Mech. Anal.
{\bf 45}
(1972),
222--240.

\bibitem{Fan_Jiang_Ni_08}
J.S. Fan, S. Jiang, G.X. Ni,
\emph{On regularity criteria for the $n$-dimensional Navier-Stokes equations in terms of the pressure,}
J. Differential Equations
{\bf 244}
(2008),
2963--2979.

\bibitem{Giga_86}
Y. Giga,
\emph{Solutions for semilinear parabolic equations in $L^p$ and regularity of weak solutions of the Navier-Stokes system,}
J. Differential Equations
{\bf 62}
(1986),
186--212.

\bibitem{HeXin_05}
C. He, Z. Xin,
\emph{On the regularity of weak solutions to the magnetohydrodynamic equations,}
J. Differential Equations
{\bf 213}
(2005),
235--254.

\bibitem{Hopf}
E. Hopf,
\emph{\"Uer die Anfangswertaufgabe f\"ur die hydrodynamischen Grundgleichungen,}
Math. Nachr.
{\bf 4}
(1951),
213-231.

\bibitem{Kim_10}
J.M. Kim,
\emph{On regularity criteria of the Navier-Stokes equations in bounded domains,}
J. Math. Phys.
{\bf 51}
(2010),
053102.

\bibitem{Kukavica_Ziane_07_one direction}
I. Kukavica, M. Ziane,
\emph{Navier-Stokes equations with regularity in one direction.}
J. Math. Phys.
{\bf 48}
(2007), 065203.

\bibitem{Kukavica_Ziane_06_one component}
I. Kukavica, M. Ziane,
\emph{One component regularity for the Navier-Stokes equations.}
Nonlinearity
{\bf 19}
(2006),
453--469.

\bibitem{Leray_34}
J. Leray,
\emph{Sur le mouvement d'un liquide visqueux emplissant l'espace.}
Acta Math.
{\bf 63}
(1934),
193--248.

\bibitem{Neustupa_Novotny_Penel_99}
J. Neustupa, A. Novotn\'y, P. Penel,
\emph{An interior regularity of a weak solution to the Navier-Stokes equations in dependence on one component of velocity,}
Topics in Mathematical Fluid Mechanics,
Quaderni di Matematica, Dept. Math., Seconda University, Napoli, Caserta,
Vol. 10, pp. 163--183 (2002);
see also
{\sl A remark to interior regularity of a suitable weak solution to the Navier-Stokes equations,}
CIM Preprint No. 25, 1999.

\bibitem{Neustupa_Penel}
J. Neustupa, P. Penel,
\emph{Anisotropic and geometric criteria for interior regularity of weak solutions to the $3$D Navier¨CStokes equations,}
in Mathematical Fluid Mechanics (Recent Results and Open Problems), Advances in Mathematical Fluid Mechanics, edited by J. Neustupa, and P. Penel (Birkh\"auser, Basel-Boston-Berlin, (2001), pp. 239--267.

\bibitem{Penel_Pokorny_10}
P. Penel, M. Pokorn\'y,
\emph{On anistropic regualarity criteria for the solutions to 3D Navier-Stokes equations,}
J. Math. Fluid Mech.
doi:10.1007/s00021-010-0038-6.

\bibitem{Penel_Pokorny}
P. Penel, M. Pokorn\'y,
\emph{Some new regularity criteria for the Navier¨CStokes equations containing the gradient of velocity,}
Appl. Math.
{\bf 49}
(2004),
483--493.

\bibitem{Prodi_59}
G. Prodi,
\emph{Un teorema di unicit\'a per le equazioni di Navier-Stokes,}
Ann. Mat. Pura Appl.
{\bf 48}
(1959),
173--182.

\bibitem{Serrin_62}
J. Serrin,
\emph{On the interior regularity of weak solutions of the Navier-Stokes equations,}
Arch. Rational Mech. Anal.
{\bf 9}
(1962),
187--191.

\bibitem{Serrin_63}
J. Serrin,
\emph{The initial value problems for the Navier-Stokes equations,}
in Nonlinear Problems, edited by R. E. Langer (University of Wisconsin Press, Madison, WI, (1963).

\bibitem{Sohr_83}
H. Sohr,
\emph{Zur Regularit\"atstheorie der instation\"aren Gleichungen von Navier-Stokes, (German) [On the regularity theory of the nonstationary Navier-Stokes equations]}
Math. Z.
{\bf 184}
(1983),
359--375.

\bibitem{Sohr_Wahl_86}
H. Sohr, W. von Wahl,
\emph{On the regularity of the pressure of weak solutions of Navier-Stokes equations,}
Arch. Math. (Basel)
{\bf 46}
(1986),
428--439.

\bibitem{Temam_77_NSE}
R. Temam,
{\sl Navier-Stokes equations, Theory and numerical analysis,}
North-Holland,
1977.

\bibitem{WuJH_03}
J.H. Wu,
\emph{Generalized MHD equations,}
J. Differential Equations
{\bf 195}
(2003),
284--312.

\bibitem{WuJH_08}
J.H. Wu,
\emph{Regularity criteria for the generalized MHD equations,}
Comm. Partial Differential Equations
{\bf 33}
(2008),
285--306.

\bibitem{Zhou_05_direction vorticity}
Y, Zhou,
\emph{A new regularity criterion for the Navier-Stokes equations in terms of the direction of vorticity,}
Monatsh. Math.
{\bf 144}
(2005),
251--257.

\bibitem{Zhou_02_grad u3}
Y. Zhou,
\emph{A new regularity criterion for the Navier-Stokes equations in terms of the gradient of one velocity component.}
Methods Appl. Anal.
{\bf 9}
(2002),
563--578.

\bibitem{Zhou_05_one component}
Y. Zhou,
\emph{A new regularity criterion for weak solutions to the Navier-Stokes equations,}
J. Math. Pures Appl.
{\bf 84}
(2005),
1496--1514.

\bibitem{Zhou_05_direction_ANZIAM}
Y. Zhou,
\emph{ Direction of vorticity and a new regularity criterion for the Navier-Stokes equations.}
ANZIAM J.
{\bf 46}
(2005),
309--316.

\bibitem{Zhou_06_grad pressure}
Y. Zhou,
\emph{On a regularity criterion in terms of the gradient of pressure for the Navier-Stokes equations in $\R^n$,}
Z. Angew. Math. Phys.
{\bf 57}
(2006),
384--392.

\bibitem{Zhou_06_pressure}
Y. Zhou,
\emph{On regularity criteria in terms of pressure for the Navier-Stokes equations in $\R^3$,}
Proc. Amer. Math. Soc.
{\bf 134}
(2006),
149--156.

\bibitem{Zhou_06_MHD pressure}
Y. Zhou,
{\sl Regularity criteria for the $3$-D MHD equations in terms of the pressure,}
Internat. J. Non-Linear Mech.
{\bf 41}
(2006),
1174--1180.

\bibitem{Zhou_04_pressure}
Y. Zhou,
\emph{Regularity criteria in terms of pressure for the $3$-D Navier-Stokes equations in a generic domain,}
Math. Ann.
{\bf 328}
(2004),
173--192.

\bibitem{Zhou_05_MHD_grad u}
Y. Zhou,
\emph{Remarks on regularities for the 3D MHD equations,} Discrete Contin. Dyn. Syst.
{\bf 12}
(2005),
881--886.

\bibitem{Zhou_Pokorny_09_grad one component}
Y. Zhou, M. Pokorn\'y,
\emph{On a regularity criterion for the Navier-Stokes equations involving gradient of one velocity component,}
J. Math. Phys.
{\bf 50}
(2009),
123514.

\bibitem{Zhou_Pokorny_09_one component}
Y. Zhou, M. Pokorn\'y,
\emph{On the regularity to the solutions of the Navier-Stokes equations via one velocity component,}
Nonlinearity
{\bf 23}
(2010),
1097--1107.

\bibitem{Zhang_Xicheng_08}
X.C. Zhang,
\emph{A regularity criterion for the solutions of $3$D Navier-Stokes equations,}
J. Math. Anal. Appl.
{\bf 346}
(2008),
336--339.

\bibitem{Zhang_Chen_05}
Z.F. Zhang, Q.L. Chen,
\emph{Regularity criterion via two components of vorticity on weak solutions to the Navier-Stokes equations in $\R^3$,}
J. Differential Equations
{\bf 216}
(2005),
470--481.

\bibitem{Zhang_10}
Z.J. Zhang,
\emph{A Serrin-type regularity criterion for the Navier-Stokes equations via one velocity component}, to appear.

\bibitem{Zhang_10_JMAA}
Z.J. Zhang,
\emph{Remarks on the regularity criteria for generalized MHD equations,}
J. Math. Anal. Appl.
{\bf 375}
(2011),
799--802.

\bibitem{Zhang_10_some}
Z.J. Zhang, M. Lu, L.D. Ni,
\emph{Some Serrin-type regularity criteria for weak solutions to the Navier-Stokes equations,}
submitted.

\bibitem{Zhang_10_two}
Z.J. Zhang,
\emph{Two new regularity criteria for the $3$D Navier-Stokes equations via two entries of the velocity gradient tensor,}
submitted.

\end{thebibliography}
\end{document}